\documentclass[11pt]{article}
\usepackage{amssymb,latexsym}
\usepackage{fullpage}
\usepackage{array}

\usepackage{epsfig}                
\def\scfig #1 #2 {\resizebox{#2}{!}{\includegraphics{#1}}}
\usepackage{fullpage}

\usepackage{xypic}
\usepackage{psfrag}

\newcommand{\lr}{\mbox{$\longrightarrow$}}

\newcommand{\be}{\begin{equation}}
\newcommand{\ee}{\end{equation}}

\newtheorem{guess}{Theorem}[section]
\newcommand{\bth}{\begin{guess}$\!\!\!${\bf .}~}
\newcommand{\eeth}{\end{guess}}
\renewcommand{\bar}{\overline}
\newtheorem{propo}[guess]{Proposition}
\newcommand{\bpropo}{\begin{propo}$\!\!\!${\bf .}~}
\newcommand{\epropo}{\end{propo}}

\newtheorem{lema}[guess]{Lemma}
\newcommand{\blem}{\begin{lema}$\!\!\!${\bf .}~}
\newcommand{\elem}{\end{lema}}

\newtheorem{defe}[guess]{Definition}
\newcommand{\bdefe}{\begin{defe}$\!\!\!${\bf .}~}
\newcommand{\edefe}{\end{defe}}

\newtheorem{coro}[guess]{Corollary}
\newcommand{\bcor}{\begin{coro}$\!\!\!${\bf .}~}
\newcommand{\ecor}{\end{coro}}

\newtheorem{rema}[guess]{Remark}
\newcommand{\brem}{\begin{rema}$\!\!\!${\bf .}~\rm}
\newcommand{\erem}{\end{rema}}

\newtheorem{exam}[guess]{Example}
\newcommand{\beg}{\begin{exam}$\!\!\!${\bf .}~\rm}
\newcommand{\eeg}{\end{exam}}

\newcommand{\bz}{\mathbb{Z}}
\newcommand{\cl}{\mathcal{L}}
\newcommand{\cj}{\mathcal{L}_+}
\newcommand{\cJ}{\bar{\mathcal{L}}_+}

\newcommand{\br}{\mathbb{R}}
\newcommand{\wt}{\widetilde}

\renewcommand{\phi}{\varphi}
\renewcommand{\hom}{$Homeo$}

\begin{document}
\title{Stone-\v{C}ech compactifications and homeomorphisms of 
products of the long line}
\author{Veerendra Vikram Awasthi and Parameswaran Sankaran  \\
Institute of Mathematical Sciences\\
CIT Campus, Chennai 600 113, INDIA \\
E-mail:
 {\tt vvawasthi@imsc.res.in}\\
{\tt sankaran@imsc.res.in}\\[2mm]}
\date{}
\maketitle
\thispagestyle{empty}

\noindent 

\noindent {\bf Abstract:} In this note we shall prove that the
Stone-\v{C}ech compactification of $\mathcal{L}^n$ is the space
$\bar{\mathcal{L}}^n$ where $\bar{\mathcal{L}}$ is the extended long line,
namely, $\mathcal{L}$ together with its ends $\pm \Omega$. We give a
similar description for the Stone-\v{C}ech compactification of 
the cartesian power of the semi-closed half-long line $\mathcal{L}_+$. As
an application we show that any torsion subgroup of the group of
all homeomorphisms of $\cj^n$ (resp. $\cl^n$) is isomorphic to a
subgroup of the symmetric group $S_n$  (resp. the
semidirect product $(\bz/2\bz)^n\ltimes S_n$).

\section{Introduction}
The well-known classification theorem for one-dimensional manifolds
is that there are exactly two non-metrizable manifolds of dimension
$1$ (besides the two metrizable ones $\br$ and $\mathbb{S}^1$)
namely the Alexandrof's long line and the half-long line. The
half-long line is described as the space $\mathcal{L}_+\setminus\{0\}$ where 
$\mathcal{L}_+:=[0,1)\times
S_\Omega$ with lexicographic order topology  and $0:=\{(0,0)\}$ is the smallest element. Here the
$S_\Omega$, the set of all ordinals which are less
than $\Omega$, is given the order topology, $\Omega$ being the smallest uncountable ordinal. 
The long line is then the space $\mathcal{L}=\mathcal{L}_-\cup \mathcal{L}_+$
glued at $\{0\}$ where $\mathcal{L}_-$ stands for $\mathcal{L}_+$
with its order reversed. The topology of these spaces have been 
studied mainly for their pathological properties in contrast with
the well-behaved metrizable manifolds. However, we 
shall see that the Stone-\v{C}ech compactifications of $\cj^n$ and $\cl^n$ 
have elegant descriptions.  Our main application is to problem of  
the characterizing,  up to isomorphism, the  torsion subgroups of the group of all homeomorphisms of $\cj^n$ 
and $\cl^n$.        

Recall that the Stone-\v{C}ech
compactification of $\mathcal{L}_+$ equals its one-point
compactification, which is the extended half-long line
$\bar{\mathcal{L}}_+:=\mathcal{L}_+\cup \{\Omega\}$, and that of 
$\mathcal{L}$ is the extended long line $\bar{\mathcal{L}}:=\cj\cup
\{\Omega, -\Omega\}$, which is also its Freudenthal's
end-compactification. 

\bth\label{main}
The Stone-\v{C}ech compactifications of $\mathcal{L}_+^n$ and
$\mathcal{L}^n$ are respectively $\bar{\mathcal{L}}_+^n$ and
$\bar{\mathcal{L}}^n$ for $n\geq 1$.
\eeth

Although this is
given as an exercise in \cite[p.255-256]{dugundji} when $n=1,2,$ we could not
find in the literature the result stated for a general $n$.

We denote by $\hom(X)$ the group of all homeomorphisms of $X$. 
The symmetric group  $S_n$ acts on $(\bz/2\bz)^n$ by permuting the factors. We denote by 
$G_n$ the semi-direct product $(\bz/2\bz)^n\ltimes S_n$.  Note that $S_n$ acts on 
$\cj^n$ and $\cl^n$ by permuting the coordinates. Also there is an obvious involution 
on $\cl$ which yields an action of $(\bz/2\bz)^n$ on $\cl^n$, $n\geq 1$. This, together with 
the action of $S_n$, defines an action of $G_n$ on $\cl^n$. 

Our next result is:
\bth \label{rigidity} There exist surjective homomorphisms
$\Phi\colon$ {\em Homeo}$(\cj^n)\lr S_n$ and $\Psi\colon$ {\em Homeo}$(\cl^n)\lr G_n$ 
which split. Furthermore, $\ker(\Phi)$ and $\ker(\Psi)$ are
torsion-free.   \eeth

An immediate corollary is the following:
\bth \label{finitesubgroups} Let $T$ be a 
subgroup of {\em Homeo}$(\cl^n)$ (resp. {\em Homeo}$(\cj^n)$ in which every element is of finite order. Then 
$T$ is isomorphic to a subgroup of $G_n$ (resp. $S_n$). \hfill $\Box$ 
\eeth

Our proofs involve mostly elementary concepts from set-topology and
make repeated use of the property that $\mathcal{L}_+$ and
$\mathcal{L}$ are sequentially compact. 
We shall also make repeated use of  
\cite[Theorem 1, \S 5.5]{mz} in \S3.  

For values of $n\leq 3$ (resp. $n\leq 2$), order of torsion elements of $\hom(\cj^n)$ (resp. $\cl^n)$) 
have been determined by Deo and Gauld \cite{dg}.  
After this paper was completed, Satya Deo pointed out to us the preprint    
\cite{bdg} in which the orders of torsion elements of $\hom(\cj^n)$ and $\hom(\cl^n)$, $n\geq 1,$ 
have been determined, but not the {\it structure} of torsion subgroups.

\section{Proof of Theorem \ref{main}} 
We use the following notations throughout.  
As usual $I$ denotes the interval $[0,1]\subset \mathbb{R}$. 
If $M$ is a manifold with boundary (which could be empty), $\partial M$ will denote the 
boundary of $M$.  If $x\in \cl$, then $-x\in \cl$ denotes the image of $x$ under the order reversing involution $\cl\lr \cl$ which switches $\cj$ and $\cl_-$ fixing $0$.
We denote by $\delta\cJ^n$ the space $\cJ^n\setminus\cj^n$. Let $x=(x_1,\cdots, x_n)\in
\cJ^n$.  $\Omega(x)$ denotes the set $\{j|x_j=\Omega\}\subset \{1,2,\cdots, n\}$ and 
$\Omega_n$ denotes the point $(\Omega,\cdots, \Omega)\in \cJ^n$.   For any set-map $h:X\lr I$, and any subset  $J\subset \br$, $h^{-1}(J)$ will have
the obvious meaning, namely, $h^{-1}(I\cap J)$.

Let $x\in \delta \cJ^n$. 
There is an obvious embedding $\gamma_x\colon \cj\lr \cj^n$  defined as $t\mapsto x(t)$ where
$x_i(t)=t$ if $i\in \Omega(x)$ and $x_i(t)=x_i$ if $i\notin \Omega(x)$. Let $f:\cj^n\lr I$ be any continuous function.  
Since $f\circ \gamma_x:\cj\lr I$ is eventually constant, it can be extended continuously to $\cJ$; its 
value at $\Omega$ will be denoted $\lim_{t\to \Omega}f(\gamma_x(t))$.
Define $\wt{f}\colon \cJ^n\lr I$ to be the set-map which extends $f$ where $\wt{f}(x):=\lim_{t\to
\Omega} f(\gamma_x(t))$ for all $x\in \delta\cJ^n$. In the case of $\cj^n$, the proof of Theorem
\ref{main} involves showing that $\wt{f}$ is continuous.  The proof will be preceded by a couple of
observations. 

\blem \label{omegan} Let $n\geq 1$ and let $h:\cJ^n\lr I$ be any continuous map. Then $h$ is constant in a
neighbourhood of $\Omega_n=(\Omega,\ldots, \Omega)$. Also any continuous map $g:\bar{\cl}^n \lr I$ is constant in a neighbourhood of $\{-\Omega, \Omega\}^n$.  \elem

\noindent {\it Proof.}  Let $h(\Omega_n)=t_0$. Choose $\lambda_k<\Omega$ such that the basic open
set $U_k=(\lambda_k,\Omega]^n$ is contained in $h^{-1}((t_0-2^{-k}, t_0+2^{-k}))$, $k\geq
1$. Let $\lambda_0$ be the least upper bound of $\{\lambda_k\mid  k\geq 1\}$. If $p\in U_{\lambda_0}$,
then $|h(p)-t_0|\leq 2^{-k}$ for every $k\geq 1$. Thus $h(p)=h(\Omega_n)$. Proof of the assertion
concerning $g$ is entirely similar and hence omitted. \hfill $\Box$

\blem\label{continuityatomegan} Let $n\geq 1$ and let $f:\cj^n\lr I$ be any continuous map.
Then $\wt{f}$ defined above is continuous at $\Omega_n\in \cJ^n$.\elem

\noindent {\it Proof.} Let $t_0=\wt{f}(\Omega_n)=\lim_{\lambda\to \Omega}f(\lambda,\ldots,\lambda)$. 
First we shall show that given any $\epsilon>0$, there exists a $\mu<\Omega$ 
such that $(\mu, \Omega)^n\subset \wt{f}^{-1}((t_0-\epsilon, t_0+\epsilon))$. 
Indeed we shall assume that there is no such $\mu$
and arrive at a contradiction. 

Choose  $\lambda_0<\Omega$ such that $f(\lambda, \ldots, \lambda)=t_0$ for
$\lambda_0<\lambda<\Omega.$ Assume that there exists an $\epsilon>0$ such that
$\wt{f}^{-1}((t_0-\epsilon, t_0+\epsilon))$ does contain $(\lambda,\Omega)^n$ for any
$\lambda<\Omega$. Choose a $p_1=(p_{1,1},\ldots, p_{1,n})\in (\lambda_0,\Omega)^n$ such that 
$|f(p_1)-t_0|\geq \epsilon$. Set $\lambda_1:=\max\{p_{1,j}\mid 1\leq j\leq n\}$ and choose
$p_2=(p_{2,1}, \ldots, p_{2,n})$ such that $\lambda_1<p_{2,j}<\Omega ~\forall j$ and
$|f(p_2)-t_0|\geq \epsilon$. Continuing thus, we obtain a sequence
$\lambda_0<\lambda_1<\cdots<\Omega$ and a sequence of points $p_1,p_2, \cdots$ in $\cj^n$
such that $\lambda_k=\max\{p_{kj}\mid 1\leq j\leq n\}, \lambda_k<p_{k+1, j}, 1\leq j\leq n,$
and $|f(p_j)-t_0|\geq \epsilon$ for all $j\geq 1$. It is evident that the sequence $(p_k)$
has  limit $p=(\nu, \ldots, \nu)\in \cj^n$ where $\nu=\lim \lambda_k$. Hence by
continuity of $f$ at $p$, we obtain $|f(p)-t_0|\geq \epsilon$. But this is a contradiction 
since, by our choice of $\lambda_0$, we have $f(\nu)=t_0$ as $\nu>\lambda_0.$ 
Hence $\wt{f}^{-1}((t_0-\epsilon,t_0+\epsilon))$ must contain $(\mu,\Omega)^n$ 
for some $\mu<\Omega$ as claimed.  Now choose such a $\mu.$  

We claim that $(\mu,\Omega]^n\subset \wt{f}^{-1}((t_0-\epsilon, t_0+\epsilon))$. 
To see this observe that if $x\in (\mu,\Omega]^n$ with $\Omega(x)\neq \emptyset$, 
then $\gamma_x(t)\in (\mu,\Omega)^n$ for all $t>\mu$. Therefore 
$|t_0-f(\gamma_x(t))|<\epsilon$ for $t>\mu$. Since $f\circ \gamma_x$ is 
eventually constant, we have $\wt{f}(x)=\lim_{t\to\Omega}f(\gamma_x(t))=
f(\gamma_x(\lambda))$ for some $\lambda>\mu$ and so $|\wt{f}(x)-t_0|<\epsilon$.
This completes
the proof. \hfill $\Box$

\noindent
{\it Proof of Theorem \ref{main}:} We prove the theorem by induction on $n$. As remarked
already, it is a well-known property when $n=1$.  We shall only consider the case of $\cj^n$. The
proof in the case of $\cl^n$ involves only obvious and routine changes.

Now suppose that $n\geq 2$ and that the theorem holds for all positive integers up to $n-1$. Let
$p=(p_1,\cdots, p_n)\in \delta \cJ^n$. Then $\Omega(p)=\{j\mid p_j=\Omega\}\neq 
\emptyset$.
Continuity at $p=\Omega_n$ having  been  established in Lemma \ref{continuityatomegan}, we assume
that $1\leq k:=n-\#\Omega(p)<n$.

\begin{center}
\scfig{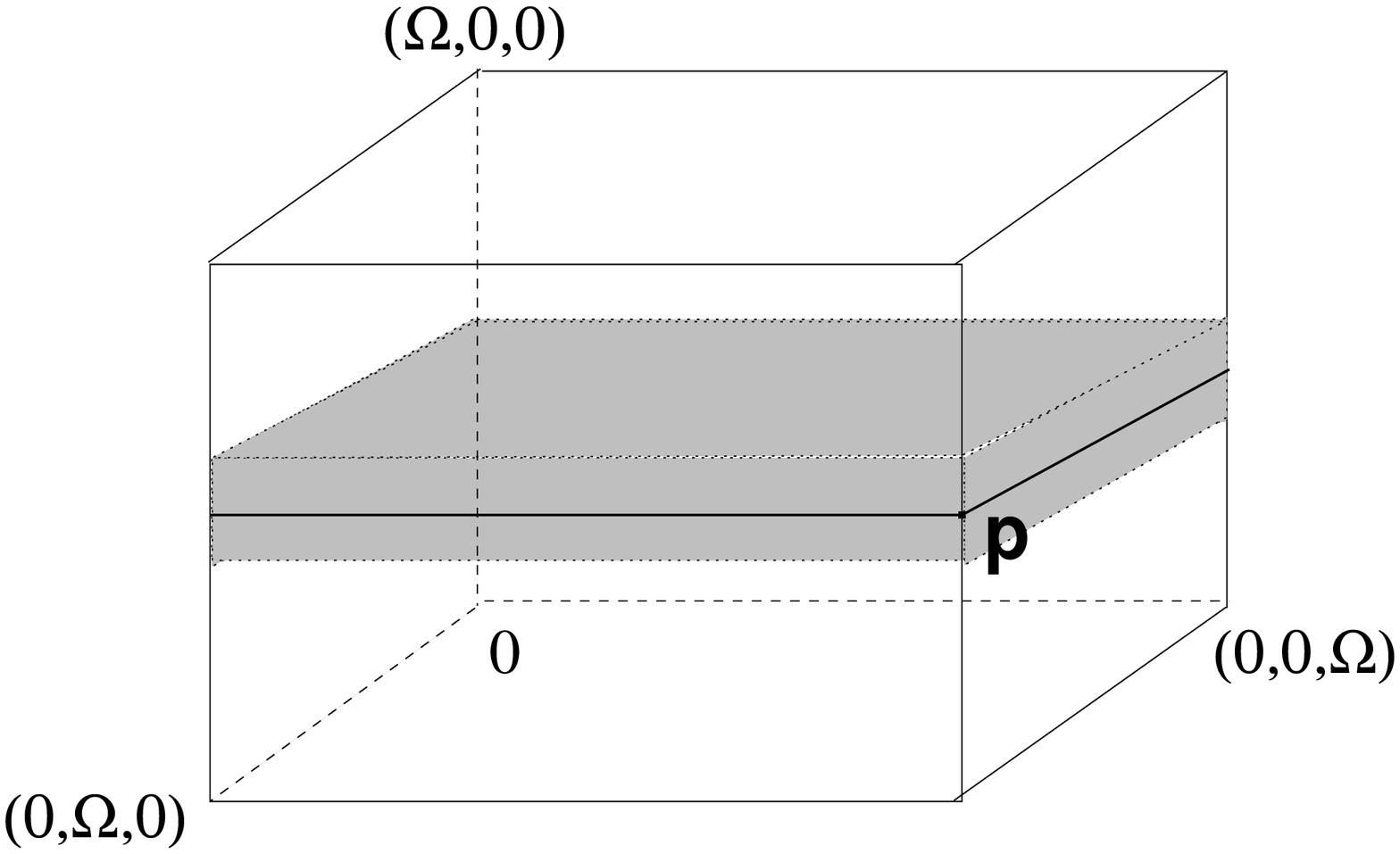} {8cm}
\end{center}
\centerline{Figure 1: $J\times \cJ^{n-k}$. }

For convenience we shall suppose that $\Omega(p)=\{j\mid k+1\leq j\leq n\}$.  Set $a:=(p_1,\ldots, 
p_k)\in \cj^k$. 
Choose a `small'  basic
open set $J\subset \cj^k$ containing $a$, i.e., $J$ is a product of bounded open intervals in $\cj$.  Let  $A$ be a countable dense subset of $J$ such that $a\in A$. For any $b\in J$, let $f_b$ denote the
restriction of $f$ to the slice $Y_b:=\{b\}\times \cj^{n-k}$. By inductive hypothesis
$\wt{f}_b:\{b\}\times \cJ^{n-k}\lr I$ is continuous. Now by  Lemma \ref{omegan} $\wt{f}_b$ is constant in a
neighbourhood of $(b,\Omega_{n-k})$. For each $b\in A$, choose an element  $\lambda_b<\Omega$ such that 
$\wt{f}_b$ is constant on $\{b\}\times U_{\lambda_b}$ where $U_{\lambda}=\{ q\in \cJ^{n-k}\mid
q_i>\lambda ~\forall i\}$. Since the $\{\lambda_b\mid b\in A\}$ is a bounded countable set, it has
a least upper bound $\mu<\Omega$. Now $\wt{f}$ restricted to $A\times U_\mu$ factors through the
projection $A\times U_\mu\lr A$.  That is, $\wt{f}(b,x)=f(b,\mu)$ for all $b\in A$ and $x\in [\mu, \Omega]$. 
Since $A\times U_\mu$ is dense in $J\times U_\mu$, this last equality holds for all $b\in J$ and $x\geq \mu$ by continuity of $f$.  It follows that $\wt{f}$ is continuous at $p$. \hfill $\Box$   

We remark that by our proof of  
Theorem \ref{main}  and Lemma \ref{omegan},  {\it any}  continuous map
$\cj^n\lr I$ is constant in a basic open set $(\mu,\Omega)^n$ for some
$\mu$. In fact, as a corollary of the proof of Theorem \ref{main}, we obtain the
following proposition.
Let $p=(\alpha,\cdots,\alpha)\in \cj^n$ be
any point on the diagonal. Then there is a retraction $r_\alpha:\cJ^n\lr
I_\alpha^n$ where $I_\alpha=[0,\alpha]^n\cong I^n$. Specifically, let
$r_\alpha$ be the product of retractions $\cJ\lr I_\alpha$ which sends
any $\beta\geq\alpha $ to $\alpha$. Similarly, one
has the retraction $\bar{\cl}\lr J_\alpha:=[-\alpha,\alpha]\cong I$
defined by $\beta\mapsto \alpha$ and
$-\beta\mapsto -\alpha$
if $\beta\geq \alpha$. This defines a retraction
$\rho_\alpha \colon \bar{\cl}^n\lr J_\alpha^n$.

\bpropo
Any continuous map $f:\cJ^n\lr I$ factors through the retraction
$r_\alpha$ for some $\alpha<\Omega$.
Similarly, any continuous map $g\colon \bar{\cl}^n\lr J_\alpha$
factors through the retraction $\rho_\alpha$ for some $\alpha<\Omega$.
\epropo

\noindent
{\it Proof.}  We shall prove the statement for $f$, the case of
$g$ being analogous.

Let $p\in \delta \cJ^n=\cJ^ n\setminus \cj^n. $ 
If $p=\Omega_n$, let $\mu\in \cj$ be such that $f|(\mu, \Omega]^n$ is constant. 

Now let $p\neq \Omega_n$ so that  $0<\#\Omega(p)<n$.
In the course of our 
proof of Theorem \ref{main}, it was observed that there is a
basic open set $U(p):=\prod U_i(p)$ containing $p$ where
$U_i(p)$ is a bounded
interval in
$\cj$ if $i\notin \Omega(p)$ and $U_i(p)=(\mu(p),\Omega]$ if $i\in
\Omega(p)$ such that $f|U(p)$ factors
through the projection $U(p)\lr
\prod_{j\notin\Omega(p)}U_j(p)\prod_{j\in \Omega(p)}\{\mu(p)\}$.

The collection $\{U(p)\}, p\in \delta \cJ^n,$ is an open covering
of the compact space $\delta \cJ^n$. Therefore there exists
finitely many points $p_1,\cdots, p_k$ such that $\{U(p_j)\}$ covers
$\delta \bar{\cl}^n$.  Now let $\mu=\max\{\mu(p_j)\mid 1\leq j\leq
k\}$. It is straightforward to verify that $f$ factors through
$r_\mu$. \hfill $\Box$

\brem\label{homotopy}
The same arguments as in the proof of Theorem \ref{main} also shows that the Stone-\v{C}ech compactification  
of $\cj^n\times I $ is $\cJ^n\times I$ and that of $\cl^n\times I $ is $\bar{\cl}^n\times I$. 
It follows readily that if $f,g:\cj^n\lr \cj^n$ are continuous, then $f$ is homotopic to 
$g$ if and only if their extensions $\wt{f}, \wt{g}:\cJ^n\lr \cJ^n$ are. An 
analogous statement holds for self-maps of $\cl^n$. 
\erem

\section{Homeomorphisms of $\cj^n$ and $\cl^n$}

Observe that $\cJ^n,\bar{\cl}^n, n\geq 1,$ are not path connected.  The path components of $\cJ^n$ are
labelled by the set $V_n:=\{0,\Omega\}^n$. More precisely, the elements of $V_n$ are in distinct
path components and every path component of $\cJ^n$ contains a point of  $V_n$. We shall denote the
path component containing $p\in V_n$ by $X_p$. Note that $X_p=\prod U_j$ where $U_j=\cj$ if
$j\notin\Omega(p)$ and $U_j=\{\Omega\}$ otherwise. In particular, $\dim X_p=n-\#\Omega(p)$.

We obtain a directed graph $\mathcal{H}_n$ (or just $\mathcal{H}$) whose vertex set is
$\mathcal{H}^0=\{X_p\mid p\in V_n\}$ and edge set $\mathcal{H}^1=\{\epsilon(p,q)\mid X_p\subset \bar{X}_q, \dim
X_p=\dim X_q-1\}$. The edge $\epsilon(p,q)$ is oriented so that it issues from $X_p$ to $X_q$. One
has a partition of $V_n=\cup_{0\leq k\leq n}V_n(k)$  where $V_n(k)=\{p\in V_n\mid \#\Omega(p)=k\}.$
The number of edges issuing from (resp. terminating at) $X_p$ equals $\#\Omega(p)$ (resp. $\dim
X_p$). 
We let $H_n$ denote the group of all automorphisms of
the directed graph $\mathcal{H}_n$.  Every element of $H_n$ fixes $X_{\Omega_n}=\{\Omega_n\}$.  
It is not hard to see that  $H_n$ is isomorphic, via restriction, to the group 
of permutations of $\{X_p\mid p\in V_n(n-1)\}\subset \mathcal{H}_n^0$. Thus 
$H_n\cong S_n$.

Let $h:\cJ^n\lr \cJ^n$ be a homeomorphism and let $h_*$ denote the induced map on the set of
path-components of $\cJ^n$.

\bpropo\label{permutations} Any homeomorphism $h\colon \cJ^n\longrightarrow\cJ^n$ induces an
isomorphism of the directed graph $\mathcal{H}_n$. The map $h\mapsto h_*$ is a surjective 
homomorphism of groups $\Phi:\hom(\cJ^n)\lr H_n$ which splits. \epropo

\noindent {\it Proof.}  Any homeomorphism of $\cJ^n$ induces an isomorphism of the set of path
components. Furthermore this defines a homomorphism from the group $\hom(\cJ^n)$ to the group of
permutations of the set of path-components of $\cJ^n$. So $h\in \hom(\cJ^n)$ induces a bijection of
the vertex set the graph $\mathcal{H}_n$. Also $h(\bar{X}_p)\subset h(X_q)$ if $\bar{X}_p \subset
X_q$. It follows that it preserves the oriented edges of $\mathcal{H}_n$. Hence $h$ induces an
isomorphism $h_*$ of $\mathcal{H}_n$.

Clearly every homeomorphism of $\cJ^n$ fixes $\Omega_n$. 
Hence $h$ maps $V_n(n-1)$ onto $V_n(n-1)$.
Any permutation $\sigma$ of $V_n(n-1)$ is evidently realizable by the homeomorphism $h$ of $\cJ^n$
given by the {\it same} permutation of the coordinates.  Therefore $\Phi$ is surjective. This also
shows that $\Phi$ splits. \hfill $\Box$

Now consider the space $\bar{\cl}^n$. Since the path components of $\bar{\cl}$ are $-\Omega, \Omega$ and
$\cl$, $\bar{\cl}^n$ has $3^n$ path components. They are labelled by $\{-\Omega, 0, \Omega\}^n$.
The element $q\in \{-\Omega, 0, \Omega\}^n$ labels the path component $X_q=\prod U_j$ where $
U_j=\{q_j\}$ if $q_j\neq 0$ and $U_j=\cl$ if $q_j=0$. Observe that $\dim X_q=\#\{j\mid \dim
q_j=0\}$.

Let $W_n=\{-\Omega, \Omega\}^n$. Each element of $W_n$ forms a path component of $\bar{\cl}^n$. Observe
that {\it any} self-homeomorphism $h$  of $\bar{\cl}^n$ preserves $W_n$ as  the other path
components are of positive dimension.

Consider the (simple) graph $\mathcal{G}_n$ whose vertices are $X_p, ~p\in W_n$. The edges of the 
graph are $e_{p,q}$ if $p$ and $q$ differ exactly in one coordinate (where they differ by sign).  The 
group of automorphisms of  $\mathcal{G}_n$ is isomorphic to the semi-direct product $(\bz/2\bz)^n\ltimes S_n:=G_n$ where 
the actions of $(\bz/2\bz)^n$ and $S_n$ are obtained from their  obvious 
respective actions  on $W_n$.

\bpropo \label{signedpermutations} Any homeomorphism $h$ of $\bar{\cl}^n$ induces an isomorphism $h_*$ of
the graph $\mathcal{G}_n$.  Furthermore, $h\mapsto h_*$ defines a surjective homomorphism of groups
$\Psi: \hom(\bar{\cl}^n)\lr G_n$ which splits. \epropo

\noindent
{\it Proof.} As observed above, $h(W_n)=W_n$ and so $h$ induces 
a bijection of the vertices of $\mathcal{G}_n$. Suppose that
$e_{p,q}$ is an edge of $\mathcal{G}_n$, say, $p_i=\Omega=-q_i,
p_j=q_j $ for $j\neq i$. Consider $X_{p,q}
=\{x\in \bar{\cl}^n \mid x_i\in \cl, x_j=p_j, j\neq i\}\cong \cl$.
Then $h(X_{p,q})$ has to be a path component of dimension
$1$ which contains $h(p), h(q)\in W_n$ in its closure. It follows
that $h(p), h(q)$ are end points of an edge of $\mathcal{G}_n$.
Therefore $h_*$ is an isomorphism of $\mathcal{G}_n$.
 
It is evident that the homeomorphisms of $\bar{\cl}^n$ which flips the 
signs of certain coordinates forms a subgroup of $\hom(\cl^n)$ isomorphic 
to $(\bz/2\bz)^n$. These, together with the homeomorphisms which 
permute the coordinates form a group isomorphic to $(\bz/2\bz)^n\ltimes S_n$.  It is evident that $\Psi$ maps 
this subgroup isomorphically onto $G_n$. Therefore $\Psi$ is split. \hfill $\Box$
 
\blem \label{infiniteorder}Let $M$ be a connected manifold with non-empty boundary $\partial
M$. Suppose that $f\colon M\lr M$ is a homeomorphism which restricts to the identity on
$\partial M$. Then either $f$ is the identity or is of infinite order. \elem

\noindent {\it Proof.} Consider the double $\hat{M}=M_0\cup_{\partial M}M_1$  of $M$,
obtained by gluing two copies $M_0, M_1$ of $M$ along the common boundary. Since $f|\partial
M$ is identity, it extends to a homeomorphism $f_0:\hat{M}\lr \hat{M}$ where $f_0$ is just
$f$ on $M_0$ and is identity on $M_1$. Since $f_0$ is identity on a non-empty open set of
$\hat{M}$, it follows that that $f_0$ is of infinite order (cf.  \cite[Theorem 1, \S 5.5]{mz}). Hence
$f$ has to be of infinite order. \hfill $\Box$

We are now ready to prove Theorem \ref{rigidity}.

\noindent {\it Proof of Theorem \ref{rigidity}:} In view of Theorem \ref{main} any
homeomorphism of $\cj^n$ lifts to a unique homeomorphism of $\cJ^n$ and similarly any
homeomorphism of $\cl^n$ lifts to a unique homeomorphism of $\bar{\cl}^n$. Indeed one has, in fact,
isomorphisms of groups $\hom(\cj^n)\cong \hom(\cJ^n)$ and $\hom(\cl^n)\cong \hom(\bar{\cl}^n)$.

By Proposition \ref{permutations} one has a surjective homomorphisms of  groups 
$\hom(\cJ^n) \lr H_n\cong S_n$,  which splits, defined as
$h\mapsto h_*$, where $H_n$ is the group of automorphisms of the directed 
graph $\mathcal{H}_n$.  
Similarly $\hom(\bar{\cl}^n)\lr G_n$ defined by $g \mapsto g_*$ is a
surjective homomorphism of groups by Proposition \ref{signedpermutations}.  As observed 
already, $H_n\cong S_n$ and $G_n\cong (\bz/2\bz)^n\ltimes S_n$.

To complete the proof, we need only to show that $\ker(\Phi)$ and $\ker(\Psi)$ are 
torsion-free.  When $n=1$ the statement is trivial to verify. Assume that $n>1$ and that the theorem is
valid for all dimensions up to $n-1$.  

Let $f$ be any element of $\hom(\cJ^n)$ of finite order such that $f_*\in H_n$ is trivial. Observe 
that $X_0=\cj^n$.  
We shall show that $f|\partial X_0$ is the identity homeomorphism. It would then follow, in view of Lemma \ref{infiniteorder}, that  
$f|X_0$ is identity and so $f$ itself is identity as $X_0$ is dense in $\cJ^n$.

\begin{center}
\scfig{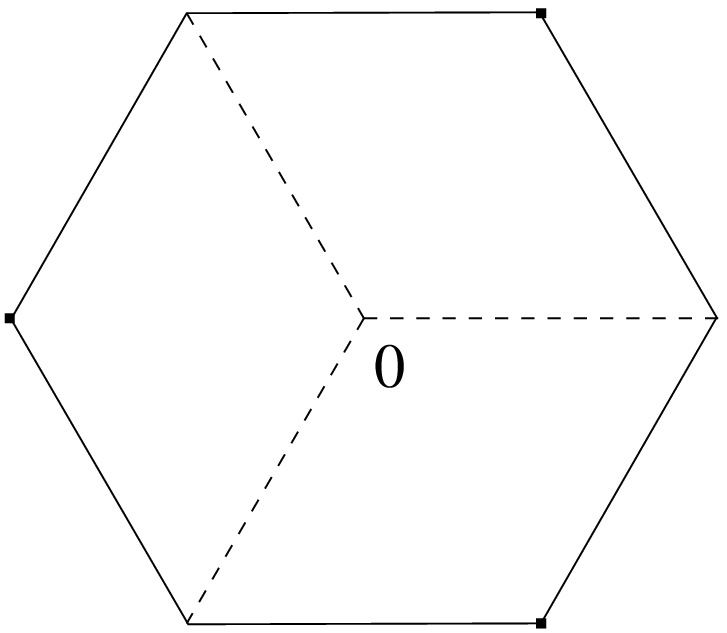} {5cm}
\end{center}
\centerline{Figure 2: $\partial X_0$.}
 
Note that 
$f$ maps each path component of $\cJ^n$ to itself and $f_p:=f|X_p$ is a finite order
homeomorphism of $X_p$ for each $p\in V_n$.  Let $k=\#\Omega(p)$. If $k>0$, then
$X_p\cong \cj^{n-k}$ and furthermore $f_p$ induces the identity map of the directed 
graph $\mathcal{H}_{n-k}$ associated to $X_p$.  Hence, by induction 
hypothesis, $f_p$ is the identity map of $X_p$. Now let $k=0$.
In this case $p=0$ and $\partial X_0$ is {\it not} homeomorphic to $\cj^n$ and so we cannot 
apply inductive hypothesis directly to conclude that $f|\partial X_0$ is the identity map. 
However, note that the points $p\in V_n(n-1)$
in the closure $\bar{\partial X}_0$ of $\partial X_0$ have basic neighbourhoods in $\bar{\partial X}_0$ 
whose closures are homeomorphic to $\cJ^{n-1}$.
Consider, say, the point $p=(\Omega, \cdots,\Omega, 0)\in \bar{\partial X}_0$. 
By what has been shown already, we have $f(p)=p$. 
Since $f$ has finite order,  
 there exists a neighbouhood $U_0\subset \bar{\partial X}_0$ of $p$ which is invariant under $f$.
Choose a $\lambda_0<\Omega$ such that the basic open set 
$B_0:=(\lambda_0,\Omega]^{n-1}\times\{0\}$ is 
contained in $U_0$. The open set $\cap_{0\leq j < r}f^r(B_0)=:U_1$ is invariant under $f$ where 
$r$ is the order of $f$. 
Repeating this argument, we get a sequence $\lambda_0<\lambda_1<\cdots <\Omega$ in $\cj$ and open 
sets $U_0\supset B_0\supset U_1 \supset B_1\supset \cdots $ in $\bar{\partial X}_0$ where 
$B_i=(\lambda_i,\Omega]^{n-1}\times \{0\}$ and  $U_i$ are invariant under $f$. Let $\mu<\Omega$ be
the limit of $(\lambda_i)$. Then  $C:=\cap_{k\geq 0}U_k = \cap_{k\geq 0} B_k = [
\mu,\Omega]^{n-1}\times\{0\}$ is invariant under $f$ and $C\cong\cJ^{n-1}$. Note that $(f|C)_*$
is the identity automorphism  of the directed graph $\mathcal{H}_{n-1}$ 
associated to $C$ since $f(x)=x$ for all $x\in X_q\cap C$ for all $q\in V_n(k), k\geq 1.$ Since $f|C$ is of finite order, by induction hypothesis we conclude that $f|C$ is 
trivial. Now by \cite[Theorem 1, \S5.5]{mz} we conclude that $f$ is the identity map.

Proof in the case of $\cl^n$ is similar and we merely give an outline.  Let $g$ be a finite order 
homeomorphism of $\bar{\cl}^n$ which induces the identity automorphism of $\mathcal{G}_n$. 
Let $X$ be any path component of $\bar{\cl}^n$ which is of dimension less than $n$. If $X$ is zero-dimensional, then it is point-wise fixed by $g$ as $g_*$ is the identity. Otherwise. 
$X\cong \cl^k, 1\leq k<n$, the map $g|
X$ is of finite order and induces the identity map of the graph associated to $X$. Hence, 
by induction hypothesis, $g|X$ is identity. Thus  $g|(\delta\bar{\cl}^n)$ is the identity.
Now consider $\Omega_n\in \bar{\cl}^n$. Proceeding as in the construction of $C$ above, we obtain 
a $\mu<\Omega$ such that $g(D)=D$ where $D:=[\mu,\Omega]^n\subset \bar{\cl}^n. $  Note that 
$D\cong \cJ^n$. Now $g|D$ is a finite order element and it induces identity map of the 
directed graph $\mathcal{H}_n$. Hence, by what has been shown already, $g|D$ is identity. 
Now by \cite[Theorem 1, \S 5.5]{mz}, $g|\cl^n$ is identity.  \hfill $\Box$

\brem
(i) Let $h:\cJ^n\lr \cJ^n$ be a homeomorphism. The induced automorphism $h_*$ of the directed graph $\mathcal{H}_n$ 
determines and is determined by the  isomorphism $H_0(h):H_0(\cJ^n;\bz)\lr H_0(\cJ^n;\bz)$ induced by $h$ in $0$-th singular homology.  Since $H_0(\cJ^n)$ is the free abelian group on the set of vertices of $\mathcal{H}_n$, 
we see that  the elements of $S_n\subset \hom(\cJ^n)$ are in distinct homotopy classes. In view of 
Remark \ref{homotopy}, we see that  distinct elements of $S_n\subset \hom(\cj^n)\cong \hom(\cJ^n)$ are in distinct homotopy classes. 
The last statement also holds for the group $G_n\subset \hom(\cl^n)$ and can be seen by a similar argument.  
See also \cite[p. 44]{bail} and \cite{bdg}.
It follows that the homomorphisms $\Phi$ and 
$\Psi$ factor through the mapping class groups 
of $\cj^n$ and $\cl^n$ respectively.  It is shown in 
\cite{bdg} that mapping class groups of $\cj^n$ and
$\cl^n$ are isomorphic to $H_n$ and $G_n$ 
respectively. Thus $\ker(\Phi)$ and  $\ker(\Psi)$ 
consist precisely of homeomorphisms which isotopic 
to respective identity maps.

(ii) It is an interesting problem to classify {\it conjugacy classes} of finite subgroups of $\hom(\cl^n)$ and of $\hom(\cj^n)$.  
\erem

\noindent 
{\bf Acknowledgements:}  The authors thank Satya Deo for a talk he gave 
at the  Ramanujan Institute for Advanced Study in Mathematics,  Chennai, in March 2008, 
which initiated our interest in homeomorphism groups of 
non-metrizable manifolds. We  
thank him also for providing us a copy of \cite{dg}.  The second author gratefully acknowledges 
financial support from Abdus Salam International Centre for Theoretical 
Physics, Trieste, Italy, during his visit in the autumn of 2008, where part of this work was 
done. 


\end{document}